\documentclass[a4paper,12pt,reqno]{amsart}

\usepackage{amsmath,amssymb,amsthm}
\usepackage{latexsym}
\usepackage{color}
\usepackage{graphicx}
\usepackage{mathrsfs}
\usepackage{enumerate}
\usepackage[abbrev]{amsrefs}
\usepackage[T1]{fontenc}
\usepackage{mathtools}
\mathtoolsset{showonlyrefs=true}

\allowdisplaybreaks[4]

\theoremstyle{plain}
\newtheorem{thm}{Theorem}[section]
\newtheorem{lemm}[thm]{Lemma}

\theoremstyle{definition}

\newtheorem{rem}[thm]{Remark}

\makeatletter

\@addtoreset{equation}{section}
\makeatother

\renewcommand{\div}{\operatorname{div}}
\newcommand{\dB}{\dot{B}}
\newcommand{\dH}{\dot{H}}

\newcommand{\supp}{\operatorname{supp}}

\renewcommand{\leq}{\leqslant}
\renewcommand{\geq}{\geqslant}

\newcommand{\pnabla}{{\nabla}^{\perp}}

\newcommand{\n}[1]{{\left\|#1\right\|}}

\newcommand{\lp}[1]{\left[#1\right]}
\newcommand{\Mp}[1]{\left\{#1\right\}}
\renewcommand{\sp}[1]{\left(#1\right)}
\newcommand{\abs}[1]{\left|#1\right|}
\newcommand{\lr}[1]{\left\langle#1\right\rangle}

\begin{document}
\title[$2$D stationary quasi-geostrophic equation]
{Stationary solutions to the critical and super-critical quasi-geostrophic equation in the scaling critical Sobolev space}
\author{Mikihiro Fujii}
\address{Graduate School of Science, Nagoya City University, Nagoya, 467-8501, Japan}
\email{fujii.mikihiro@nsc.nagoya-cu.ac.jp}
\begin{abstract}
We consider the stationary problem for the quasi-geostrophic equation with the critical and super-critical dissipation and prove the unique existence of small solutions for given small external force in the scaling critical Sobolev spaces framework.
Moreover, we also show that the data-to-solution map is continuous.
Since the critical and super-critical case involves the derivative loss, which affects the class of the continuity of the data-to-solution map, we reveal that the map is no longer uniform continuous, in contrast to the sub-critical case, where the Lipschitz continuity holds.
\end{abstract}
\maketitle


\section{Introduction}\label{sec:intro}
Let us consider the stationary problem of the quasi-geostrophic equation:
\begin{align}\label{eq:s-QG}
    \begin{cases}
        (-\Delta)^{\alpha} \theta + v \cdot \nabla \theta = f, & \quad x \in \mathbb{R}^2, \\
        v = \nabla^{\perp}(-\Delta)^{-1/2} \theta, & \quad x \in \mathbb{R}^2,
    \end{cases}
\end{align}
where $\alpha$ is a positive constant and $\theta=\theta(x)$ is the unknown potential temperature, whereas $f=f(x)$ is the given external force.
The fractional Laplacian is defined by $(-\Delta)^{s}:=\mathscr{F}^{-1}|\xi|^{2s}\mathscr{F}$ for $s \in \mathbb{R}$, and we denote $\nabla^{\perp}:=(\partial_{x_2},-\partial_{x_1})$.
The dissipation of the equation \eqref{eq:s-QG} is given by the fractional Laplacian $(-\Delta)^{\alpha}\theta$ and the cases $\alpha>1/2$, $\alpha=1/2$, and $0<\alpha\leq 1/2$ is so called sub-critical, critical, and super-critical cases, respectively.
It is known that critical exponent $\alpha=1/2$ is the threshold that distinguishes whether equation \eqref{eq:s-QG} is semilinear or quasilinear, with sub-critical cases being semilinear and critical and super-critical cases being quasilinear.
The equation \eqref{eq:s-QG} has an invariant scaling structure;
if $\theta$ is a solution to \eqref{eq:s-QG} with some external force $f$, then the scaled function $\theta_{\lambda}(x)=\lambda^{2\alpha - 1}\theta(\lambda x)$ also solves \eqref{eq:s-QG} with $f$ replaced by $f_{\lambda}(x):=\lambda^{4\alpha-1}f(\lambda x)$.
We say that the function spaces $D$ and $S$ for the external force and solutions, respectively, are scaling critical if $\n{f_{\lambda}}_D=\n{f}_D$ and $\n{\theta_{\lambda}}_S=\n{\theta}_S$ for all $\lambda>0$.
For example, $(D,S)=(\dH^{2-4\alpha}(\mathbb{R}^2),\dH^{2-2\alpha}(\mathbb{R}^2))$ or $(D,S)=(\dB_{p,q}^{2/p+1-4\alpha}(\mathbb{R}^2),\dB_{p,q}^{2/p+1-2\alpha}(\mathbb{R}^2))$ ($1\leq p,q \leq \infty$) are scaling critical spaces.
In the recent work of Kozono--Kunstmann--Shimizu \cite{KKS}, the solvability of the sub-critical quasi-geostrophic equation \eqref{eq:s-QG} was proved in the scaling critical Besov spaces framework. 
More precisely, they showed that for any $1/2 < \alpha < 1$, $1 \leq p < 2/(2\alpha-1)$, and $1 \leq q \leq \infty$, if the given external force $f \in \dB_{p,q}^{2/p+1-4\alpha}(\mathbb{R}^2)$ is small, there exists a unique small solution $\theta \in \dB_{p,q}^{2/p+1-2\alpha}(\mathbb{R}^2)$. Note that their result immediately implies the data-to-solution map $f \mapsto \theta$ is Lipschitz continuous.
Here, the sub-critical assumption $\alpha>1/2$ is caused by the controllability of the first order spatial derivative in the nonlinear term by the inverse fractional Laplacian $(-\Delta)^{-\alpha}$. In other word, to obtain the key bilinear estimate
\begin{align}\label{key:KKS}
    \n{(-\Delta)^{-\alpha}\sp{v_1 \cdot \nabla \theta_2}}_{\dB_{p,q}^{2/p+1-2\alpha}}
    \leq
    C
    \n{\theta_1}_{\dB_{p,q}^{2/p+1-2\alpha}}
    \n{\theta_2}_{\dB_{p,q}^{2/p+1-2\alpha}}
\end{align}
with $v_1:=\nabla^{\perp}(-\Delta)^{-1/2}\theta_1$,
the sub-critical condition $\alpha>1/2$ is required.
The aim of this paper is to investigate the solvability of \eqref{eq:s-QG} in the case of $\alpha \leq 1/2$ where \eqref{key:KKS} is difficult to prove due to the weakness of the dissipation.
Our first outcome of this article is that every small external force establishes a unique small solution to \eqref{eq:s-QG} in the scaling critical $L^2$-Sobolev space framework.
Moreover, we also show that the data-to-solution map is continuous; for any $\{f_n\}_{n \in \bar{\mathbb{N}}}$ with $\mathbb{N}:=\mathbb{N} \cup \{\infty\}$ and $f_n \to f_{\infty}$ as $n \to \infty$, the corresponding sequence $\{ \theta_n \}_{n \in \bar{\mathbb{N}}}$ of solutions satisfies $\theta_n \to \theta$ as $n \to \infty$.
Unlike the sub-critical case $\alpha>1/2$, the derivative loss due to the nonlinear term appears in the case of $\alpha \leq 1/2$ and this makes it difficult to prove the Lipschitz continuity of the data-to-solution map.
For the class of the continuity for the solution map, we show that the solution map is {\it not uniform continuous} in the scaling critical space.

Before stating our main results, we briefly recall the existing studies related to our work.
For the initial value problem of the quasi-geostrophic equation 
\begin{align}\label{eq:n-QG}
    \begin{cases}
    \partial_t\theta + (-\Delta)^{\alpha} \theta + v \cdot \nabla \theta = f, & \quad t\geq 0, x \in \mathbb{R}^2, \\
    v = \nabla^{\perp}(-\Delta)^{-1/2} \theta, & \quad t\geq 0, x \in \mathbb{R}^2,\\
    \theta(0,x) = \theta_0(x), & \quad x \in \mathbb{R}^2,
    \end{cases}
\end{align}
the local and global well-posedness are investigated by many researchers; see \cites{CC-04,CW-99,Wu2,Wu} for the sub-critical case $\alpha>1/2$, \cites{Cha-Lee-03, CCW-01,DD-08, Kis-Naz-Vol-07,Zha-07} for the critical case $\alpha=1/2$, and \cites{Bae,Cha-Lee-03, Miu-06, Che-Mia-Zha-07} for the super-critical case $0< \alpha <1/2$.
Here, we remark that for the sub-critical and critical case $\alpha \geq 1/2$, \cites{CW-99,Kis-Naz-Vol-07} proved the global well-posedness for arbitrary large initial data in the scaling critical framework, whereas the global regularity for the super-critical case remains open.
Next, we mention the known results for the stationary quasi-geostrophic equation \eqref{eq:s-QG}, on which there seem to be fewer results in comparison with the initial value problem on the non-stationary quasi-geostrophic equation \eqref{eq:n-QG}.
Dai \cite{Dai-15} showed the existence of weak solution $\theta \in H^{\frac{1}{2}}(\mathbb{R}^2)$ for small data $f \in W^{\frac{1}{2},4}(\mathbb{R}^2) \cap L^{\infty}(\mathbb{R}^2)$ satisfying $\widehat{f}(\xi)=0$ around $\xi = 0$. 
Moreover, she removed the condition $\widehat{f}(\xi)=0$ around $\xi = 0$ in the case of $1/2 < \alpha < 1$ and proved the existence of weak solutions $\theta \in H^{\alpha}(\mathbb{R}^2)$ for small data $f \in W^{1-\alpha,4}(\mathbb{R}^2) \cap L^{\frac{4}{2\alpha-1}}(\mathbb{R}^2)$.
For the unique existence of strong solutions in the sub-critical case $1/2 < \alpha < 1$,
Hadadifardo--Stefanov \cite{HS-21} proved that for any small external force $f \in L^{\frac{2}{2\alpha-1}}(\mathbb{R}^2)$, there exists a unique small solution $\theta \in \dot{W}^{-2\alpha,\frac{2}{2\alpha-1}}(\mathbb{R}^2)$; note that their framework is scaling critical.
Recently, Kozono--Kunstmann--Shimizu \cite{KKS} generalized the result of \cite{HS-21} and showed that for any $1/2 < \alpha < 1$, $1 \leq p < 2/(2\alpha-1)$, and $1 \leq q \leq \infty$, if the given external force $f \in \dB_{p,q}^{2/p+1-4\alpha}(\mathbb{R}^2)$ is small, there exists a unique solution $\theta \in \dB_{p,q}^{2/p+1-2\alpha}(\mathbb{R}^2)$.
However, due to the difficulty of the derivative loss in the nonlinear term, there appear to be no results on the unique existence of the solutions to \eqref{eq:s-QG} in the critical super-critical case $0<\alpha \leq 1/2$.
The aim of this paper is to provide the existence of a unique solution to \eqref{eq:s-QG} in the critical and super-critical case $0<\alpha \leq 1/2$ and investigate the continuity of data-to-solution map.

Now, we state our main result of this paper precisely.
The first result reads as follows.
\begin{thm}[Unique existence of small solutions]\label{thm:1}
    Let $0 < \alpha \leq 1/2$.
    Then, there exists a positive constant $\delta=\delta(\alpha)$ and $\eta=\eta(\alpha)$ such that for any external force $f \in \dH^{-\alpha}(\mathbb{R}^2) \cap \dH^{2-4\alpha}(\mathbb{R}^2)$ and $\n{f}_{\dH^{2-4\alpha}} \leq \delta$, \eqref{eq:s-QG} possesses a unique solution $\theta \in \dH^{\alpha}(\mathbb{R}^2) \cap \dH^{2-2\alpha}(\mathbb{R}^2)$ satisfying $\n{\theta}_{\dH^{2-2\alpha}} \leq \eta$.
\end{thm}
\begin{rem}
    Let us make some remarks on Theorem \ref{thm:1}.
    The scaling critical $L^2$-Sobolev space for $(f,\theta)$ is given by $(\dH^{2-4\alpha}(\mathbb{R}^2),\dH^{2-2\alpha}(\mathbb{R}^2))$. 
    However, as $2-2\alpha>1$, $\dH^{2-2\alpha}(\mathbb{R}^2)$ is not included in $\mathscr{S}'(\mathbb{R}^2)$; see Section \ref{sec:pre} for details.
    Therefore, we assume the auxiliary low-regularity $f \in \dH^{-\alpha}(\mathbb{R}^2)$ and $\theta \in \dH^{\alpha}(\mathbb{R}^2)$. 
    The reason of the choice of low-regularities comes from the calculation
    \begin{align}
        \n{(-\Delta)^{\frac{\alpha}{2}}\theta}_{L^2}^2 = \lr{f,\theta}_{L^2} \leq \n{f}_{\dH^{-\alpha}}\n{\theta}_{\dH^{\alpha}},
    \end{align}
    which is implied by the standard $L^2$-energy argument via the nonlinear cancellation $\lr{u \cdot \nabla \theta,\theta}_{L^2}=0$.
\end{rem}
From the point of view of the series of known results \cites{Kan-Koz-Shi-19,Tsu-19-ARMA,Tsu-19-JMAA,Tan-Tsu-Zha-pre,Li-Yu-Zhu,Fuj-24,Fuj-pre} for the stationary Navier--Stokes equations, investigating the continuity of the data-to-solution map is an important topic.
In addition, unlike the sub-critical case $1/2<\alpha<1$, where the contraction mapping argument is directly available, it is not obvious whether the data-to-solution map is continuous or not.
For the continuous dependence of the solution with respect to the external force, the following assertion holds.
\begin{thm}[Continuity of the solution map]\label{thm:2}
    Let $0 < \alpha \leq 1/2$.
    Then, the data-to-solution map $S: f \mapsto \theta$ is continuous from $\dH^{-\alpha}(\mathbb{R}^2) \cap \dH^{2-4\alpha}(\mathbb{R}^2)$ to $\dH^{\alpha}(\mathbb{R}^2) \cap \dH^{2-2\alpha}(\mathbb{R}^2)$.
    More precisely, there exists a positive constant $\delta=\delta(\alpha)$ such that for any $\{ f_n \}_{n \in \bar{\mathbb{N}}} \subset \dH^{-\alpha}(\mathbb{R}^2) \cap \dH^{2-4\alpha}(\mathbb{R}^2)$ with $\bar{\mathbb{N}}:= \mathbb{N} \cup \{ \infty \}$, satisfying $\n{f_n}_{\dH^{2-4\alpha}} \leq \delta$ for all $n \in \bar{\mathbb{N}}$ and 
    \begin{align}
        \lim_{n \to \infty} 
        \n{f_n - f_{\infty}}_{\dH^{-\alpha} \cap \dH^{2-4\alpha}} = 0,
    \end{align}
    the corresponding solution sequence $\{ \theta_n \}_{n \in \bar{\mathbb{N}}} \subset \dH^{\alpha}(\mathbb{R}^2) \cap \dH^{2-2\alpha}(\mathbb{R}^2)$ constructed in the Theorem \ref{thm:1} satisfies 
    \begin{gather}
        \n{\theta_n-\theta_{\infty}}_{\dH^{\alpha}}
        \leq
        C
        \n{f_n - f_{\infty}}_{\dH^{-\alpha}}, \label{Lip}\\
        \lim_{n \to \infty}
        \n{\theta_n-\theta_{\infty}}_{\dH^{2-2\alpha}}
        =0 \label{not-Lip}
    \end{gather}
    with some positive constant $C$ depending only on $\alpha$.
\end{thm}
From the result of Kozono--Kunstmann--Shimizu \cite{KKS}, we immediately see that the solution map $S$ is Lipschitz continuous from $\dH^{2-4\alpha}(\mathbb{R}^2)$ to $\dH^{2-2\alpha}(\mathbb{R}^2)$ for the subcritical case $1/2<\alpha<1$.
In comparison to this, 
we may be slightly generalize \eqref{Lip}; for sufficiently small $f,g \in \dH^{-\alpha}(\mathbb{R}^2) \cap \dH^{2-4\alpha}(\mathbb{R}^2)$, the corresponding small solutions $\theta[f],\theta[g] \in \dH^{\alpha}(\mathbb{R}^2) \cap \dH^{2-2\alpha}(\mathbb{R}^2)$ satisfies
\begin{align}
    \n{\theta[f] - \theta[g]}_{\dH^{\alpha}}
    \leq
    C\n{f-g}_{\dH^{-\alpha}}.
\end{align}
Thus, the data-to-solution map $S$ is Lipschitz continuous in $\dH^{\alpha}(\mathbb{R}^2)$.
On the other hand, there is no such information on \eqref{not-Lip} for the critical regularity $\dH^{2-2\alpha}(\mathbb{R}^2)$.
For this direction, we reveal that the solution map is no longer uniform continuous in $\dH^{2-2\alpha}(\mathbb{R}^2)$.
More precisely, we obtain the following theorem.
\begin{thm}[Non-uniform continuity of the solution map]\label{thm:3}
    For $0<\alpha < 1/2$, the data-to-solution map in Theorem \ref{thm:2} is not uniform continuous.
    More precisely, there exist two sequences $\{ f_n \}_{n \in \mathbb{N}}$ and $\{ g_n \}_{n \in \mathbb{N}}$ in $\mathscr{S}(\mathbb{R}^2)$ such that 
    \begin{align}
        \lim_{n \to \infty}
        \n{f_n - g_n}_{\dH^{-\alpha} \cap \dH^{2-4\alpha}} = 0,
    \end{align}
    the corresponding sequence $\{ \theta[f_n] \}_{n\in\mathbb{N}}$ and $\{ \theta[g_n] \}_{n\in\mathbb{N}}$ of the solutions in $\dH^{\alpha}(\mathbb{R}^2) \cap \dH^{2-2\alpha}(\mathbb{R}^2)$ satisfies 
    \begin{align}
        \liminf_{n \to \infty}
        \n{ \theta[f_n] - \theta[g_n] }_{\dH^{2-2\alpha}}
        >
        0.
    \end{align}
\end{thm}
\begin{rem}
    For the known results of the non-uniform continuity of the data-to-solution map, we refer to \cite{BL} for the Euler equations case. 
    See also \cite{ZZ} for the anisotropic Navier--Stokes equations case.
\end{rem}
In Theorem \ref{thm:3}, the critical case $\alpha=1/2$ is removed.
For this case, using the Besov spaces with the interpolation index $q=1$, we obtain the Lipschitz continuity of the data-to-solution map; indeed we have the following result.
\begin{thm}[Well-posedness for the critical case]\label{thm:4}
    Let $\alpha=1/2$ and $1 \leq p < \infty$.
    Then, there exists a positive constant $\delta=\delta(p)$ and $\eta=\eta(p)$ such that if $f \in \dB_{p,1}^{2/p-1}(\mathbb{R}^2)$ satisfies $\n{f}_{\dB_{p,1}^{2/p-1}} \leq \delta$, then \eqref{eq:s-QG} with $\alpha=1/2$ admits a unique solution $\theta \in \dB_{p,1}^{2/p}(\mathbb{R}^2)$ with $\n{\theta}_{\dB_{p,1}^{2/p}}\leq \eta$. 
    Moreover, the data-to-solution map $S:f\mapsto \theta$ is Lipschitz continuous from $\dB_{p,1}^{2/p-1}(\mathbb{R}^2)$ to $\dB_{p,1}^{2/p}(\mathbb{R}^2)$.
\end{thm}
\begin{rem}
    As Theorem \ref{thm:4} treats the case $\alpha=1/2$ in the critical Besov spaces, this corresponds to one of the end-point cases of \cite{KKS}.
\end{rem}
The proof of Theorem \ref{thm:4} is immediately follows from 
the standard contraction mapping principle via the estimate
\begin{align}
    \n{(-\Delta)^{-1/2}(v_1\cdot \nabla \theta_2)}_{\dB_{p,1}^{2/p}}
    ={}&
    \n{(-\Delta)^{-1/2}\div (v_1  \theta_2)}_{\dB_{p,1}^{2/p}}\\
    \leq{}&
    C
    \n{v_1 \theta_2}_{\dB_{p,1}^{2/p}}\\
    \leq{}&
    C
    \n{v_1}_{\dB_{p,1}^{2/p}}
    \n{\theta_2}_{\dB_{p,1}^{2/p}}\\
    \leq{}&
    C
    \n{\theta_1}_{\dB_{p,1}^{2/p}}
    \n{\theta_2}_{\dB_{p,1}^{2/p}}
\end{align}
for $\theta_1, \theta_2 \in \dB_{p,1}^{2/p}(\mathbb{R}^2)$ with $v_1:=\pnabla (-\Delta)^{-1/2}\theta_1$ and $1 \leq p<\infty$.

This paper is organized as follows.
In Section \ref{sec:pre}, we prepare the definition of function spaces and basic properties of bilinear estimates.
In Sections \ref{sec:pf_thm:1}, \ref{sec:pf_thm:2}, and \ref{sec:pf_thm:3}, we provide the proofs of Theorems \ref{thm:1}, \ref{thm:2}, and \ref{thm:3}, respectively.
\section{Preliminaries}\label{sec:pre}
In the present section, we prepare some notation and key lemmas used in this manuscript.
For $1 \leq p \leq \infty$, $L^p(\mathbb{R}^2)$ represents the usual Lebesgue space.
For $s \in \mathbb{R}$, we define the homogeneous Sobolev space of order $s$ by 
\begin{align}
    \dH^{s}(\mathbb{R}^2)
    :=
    \Mp{
    f \in \mathscr{S}'(\mathbb{R}^2)/\mathscr{P}(\mathbb{R}^2)\ ;\ 
    \n{f}_{\dH^s}:=\n{(-\Delta)^{\frac{s}{2}}f}_{L^2}<\infty
    },
\end{align}
where $\mathscr{P}(\mathbb{R}^2)$ is the set of all polynomials on $\mathbb{R}^2$.
It is well-known that for $s>1$, it holds
\begin{align}
    \dH^s(\mathbb{R}^2)
    \sim
    \Mp{f \in \mathscr{S}'(\mathbb{R}^2)\ ;\ \n{f}_{\dH^s}<\infty,\quad \sum_{j \in \mathbb{Z}} \phi_j*f=f\ {\rm in}\ \mathscr{S}'(\mathbb{R}^2)},
\end{align}
where $\{\phi_j\}_{j \in \mathbb{Z}}$ denotes the Littlewood--Paley frequency cut-off functions.
We recall the well-known product estimates in homogeneous Sobolev spaces and provide a commutator estimate.
\begin{lemm}
    The following estimates hold true:
    \begin{enumerate}
        \item 
        Let $s,s_1,s_2,s_3,s_4 \in \mathbb{R}$ satisfy $s_1<1$, $s_4<1$, and $s=s_1+s_2=s_3+s_4>0$.
        Then, there exists a positive constant $C=C(s,s_1,s_2,s_3,s_4)$ such that 
        \begin{align}
            \n{fg}_{\dH^{s-1}}
            \leq 
            C
            \n{f}_{\dH^{s_1}}
            \n{g}_{\dH^{s_2}}
            +
            C
            \n{f}_{\dH^{s_3}}
            \n{g}_{\dH^{s_4}}
        \end{align}
        for all $f \in \dH^{s_1}(\mathbb{R}^2) \cap \dH^{s_3}(\mathbb{R}^2)$ and $g \in \dH^{s_2}(\mathbb{R}^2) \cap \dH^{s_4}(\mathbb{R}^2)$.
        \item 
        Let $s_1,s_2,s_3,s_4,s_5,s_6 \in \mathbb{R}$ satisfy $s_2>0$, $s_3<2$, $s_6<1$, and $s_1+s_2=s_3+s_4=s_5+s_6>0$.
        Then, there exists a positive constant $C=C(s_1,s_2,s_3,s_4,s_5,s_6)$ such that 
        \begin{align}
            \n{\lp{(-\Delta)^{\frac{s_1}{2}}, f}g}_{\dH^{s_2-1}}
            \leq 
            C
            \n{f}_{\dH^{s_3}}
            \n{g}_{\dH^{s_4}}
            +
            C
            \n{f}_{\dH^{s_5}}
            \n{g}_{\dH^{s_6}}
        \end{align}
        for all $f \in \dH^{s_3}(\mathbb{R}^2) \cap \dH^{s_5}(\mathbb{R}^2)$ and $g \in \dH^{s_4}(\mathbb{R}^2) \cap \dH^{s_6}(\mathbb{R}^2)$, where we have set $\lp{(-\Delta)^{\frac{s_1}{2}},f}g:=(-\Delta)^{\frac{s_1}{2}}(fg)-f(-\Delta)^{\frac{s_1}{2}}g$.
    \end{enumerate}
\end{lemm}
\begin{proof}
    The first assertion is proved by the standard para-product estimates; see \cite{Bah-Che-Dan-11} for the detail.
    Therefore, we focus on the second one.
    Let $s:=s_1+s_2=s_3+s_4=s_5+s_6$.
    Let us decompose the Fourier transform of the commutator as 
    \begin{align}
        &
        \mathscr{F}\lp{\lp{(-\Delta)^{\frac{s_1}{2}}, f}g}(\xi)
        ={}
        \int_{\mathbb{R}^2}
        \widehat{f}(\eta)
        \widehat{g}(\xi-\eta)
        \sp{|\xi|^{s_1}-|\xi-\eta|^{s_1}} 
        d\eta\\
        &\quad={}
        \int_{|\eta| \leq \frac{1}{2}|\xi-\eta|}
        \widehat{f}(\eta)
        \widehat{g}(\xi-\eta)
        \int_0^1
        s_1
        |\xi-\theta\eta|^{s_1-2}\Mp{(\xi-\theta\eta)\cdot\eta}d\theta 
        d\eta\\
        &\qquad
        +
        |\xi|^{s_1} 
        \int_{\frac{1}{2}|\xi-\eta| \leq |\eta| \leq 2|\xi-\eta|}
        \widehat{f}(\eta)
        \widehat{g}(\xi-\eta)
        d\eta
        \\
        &\qquad
        +
        |\xi|^{s_1} 
        \int_{|\eta| \geq 2|\xi-\eta|}
        \widehat{f}(\eta)
        \widehat{g}(\xi-\eta)
        d\eta
        \\
        &\qquad
        -
        \int_{\frac{1}{2}|\xi-\eta| \leq |\eta| \leq 2|\xi-\eta|}
        \widehat{f}(\eta)
        \widehat{g}(\xi-\eta)
        |\xi-\eta|^{s_1}
        d\eta
        \\
        &\qquad
        -
        \int_{|\eta| \geq 2|\xi-\eta|}
        \widehat{f}(\eta)
        \widehat{g}(\xi-\eta)
        |\xi-\eta|^{s_1}
        d\eta
        \\
        &\quad=:I_1(\xi)+I_2(\xi)+I_3(\xi)+I_4(\xi)+I_5(\xi).
    \end{align}
    By the direct calculation, we see by $(1/2)|\xi| \leq |\xi-\theta \eta| \leq (3/2)|\xi|$ for $0\leq \theta \leq 1$ and $|\eta| \leq (1/2)|\xi-\eta|$ that  
    \begin{align}
        &
        |\xi|^{s_2-1}
        \abs{I_1(\xi)}
        \\
        &\quad\leq
        s_1
        \int_{|\eta| \leq \frac{1}{2}|\xi-\eta|}
        \abs{\widehat{f}(\eta)}
        \abs{\widehat{g}(\xi-\eta)}
        \int_0^1
        |\eta|
        |\xi|^{s_2-1}|\xi-\theta\eta|^{s_1-1}
        d\theta 
        d\eta\\
        &\quad\leq
        C
        \int_{|\eta| \leq \frac{1}{2}|\xi-\eta|}
        \abs{\widehat{f}(\eta)}
        \abs{\widehat{g}(\xi-\eta)}
        |\eta|
        |\xi|^{s-2}
        d\eta\\
        &\quad
        \leq
        C
        \n{f}_{\dH^{s_3}}
        \Mp{
        \int_{|\eta| \leq \frac{1}{2}|\xi-\eta|}
        \sp{|\xi-\eta|^{s-2}
        \abs{\widehat{g}(\xi-\eta)}
        |\eta|^{1-s_3}
        }^2
        d\eta}^{\frac{1}{2}},
    \end{align}
    which implies 
    \begin{align}
        &
        \n{|\xi|^{s_3}I_1(\xi)}_{L^2_{\xi}}^2
        \\
        &\quad\leq
        C
        \n{f}_{\dH^{s_3}}^2
        \int_{\mathbb{R}^2_{\eta}}
        \sp{\int_{|\eta| \leq \frac{1}{2}|\xi-\eta|}
        \sp{|\xi-\eta|^{s-2}
        \abs{\widehat{g}(\xi-\eta)}
        }^2
        d\xi
        }
        |\eta|^{2(1-s_3)}
        d\eta\\
        &\quad=
        C
        \n{f}_{\dH^{s_3}}^2
        \int_{\mathbb{R}^2_{\eta}}
        \sp{\int_{|\eta| \leq \frac{1}{2}|\xi|}
        \sp{|\xi|^{s-2}
        \abs{\widehat{g}(\xi)}
        }^2
        d\xi
        }
        |\eta|^{2(1-s_3)}
        d\eta\\
        &\quad=
        C
        \n{f}_{\dH^{s_3}}^2
        \int_{\mathbb{R}^2_{\xi}}
        \sp{\int_{|\eta| \leq \frac{1}{2}|\xi|}
        |\eta|^{2(1-s_3)}
        d\eta
        }
        \sp{|\xi|^{s-2}
        \abs{\widehat{g}(\xi)}
        }^2
        d\xi\\
        &\quad\leq
        C
        \n{f}_{\dH^{s_3}}^2
        \int_{\mathbb{R}^2_{\xi}}
        |\xi|^{2(2-s_3)}
        \sp{|\xi|^{s-2}
        \abs{\widehat{g}(\xi)}
        }^2
        d\xi\\
        &\quad=
        C
        \n{f}_{\dH^{s_3}}^2
        \int_{\mathbb{R}^2_{\xi}}
        \sp{|\xi|^{s_4}
        \abs{\widehat{g}(\xi)}
        }^2
        d\xi
        =
        C
        \n{f}_{\dH^{s_3}}^2
        \n{g}_{\dH^{s_4}}^2.
    \end{align}
    Next, we consider the estimate of $I_2(\xi)$.
    Since $|\xi|\leq3|\eta|$ if $(1/2)|\xi-\eta| \leq |\eta| \leq 2|\xi-\eta|$, we have 
    \begin{align}
        &|\xi|^{s_2-1}
        \abs{I_2(\xi)}\\
        &\quad
        \leq{}
        C|\xi|^{s-1}
        \int_{\frac{1}{2}|\eta| \leq |\xi-\eta| \leq 2|\eta|}
        |\xi-\eta|^{-s}
        |\eta|^{s_5}
        \abs{\widehat{f}(\eta)}
        |\xi-\eta|^{s_6}
        \abs{\widehat{g}(\xi-\eta)}
        d\eta
    \end{align}
    which implies 
    \begin{align}
        &
        \n{|\xi|^{s_2-1}
        \abs{I_2(\xi)}}_{L^2_{\xi}}\\
        & \quad
        \leq
        C
        \left\{
        \int_{\mathbb{R}^2_{\xi}}
        |\xi|^{2(s-1)}
        \left(
        \int_{\frac{1}{2}|\eta| \leq |\xi-\eta| \leq 2|\eta|}
        |\xi-\eta|^{-s}
        \right.\right.
        \\
        &\qquad \qquad 
        \left.\left.
        \times
        |\eta|^{s_5}
        \abs{\widehat{f}(\eta)}
        |\xi-\eta|^{s_6}
        \abs{\widehat{g}(\xi-\eta)}
        d\eta\right)^2
        d\xi
        \right\}^{\frac{1}{2}}\\
        & \quad
        \leq
        C
        \int_{\mathbb{R}^2_{\eta}}
        |\eta|^{s_5}
        \abs{\widehat{f}(\eta)}\\
        &\qquad \quad
        \times
        \Mp{
        \int_{\frac{1}{2}|\eta| \leq |\xi-\eta| \leq 2|\eta|}
        \sp{|\xi|^{s-1}|\xi-\eta|^{-s}
        |\xi-\eta|^{s_6}
        \abs{\widehat{g}(\xi-\eta)}}^2
        d\eta}^{\frac{1}{2}}
        d\xi \\
        & \quad
        \leq
        C
        \n{f}_{\dH^{s_5}}
        \Mp{
        \int_{\mathbb{R}^2_{\eta}}
        \int_{\frac{1}{2}|\eta| \leq |\xi-\eta| \leq 2|\eta|}
        \sp{|\xi|^{s-1}|\xi-\eta|^{-s}
        |\xi-\eta|^{s_6}
        \abs{\widehat{g}(\xi-\eta)}}^2
        d\eta
        d\xi}^{\frac{1}{2}}\\
        & \quad
        =
        C
        \n{f}_{\dH^{s_5}}
        \Mp{
        \int_{\mathbb{R}^2_{\xi}}
        \sp{|\xi|^{s_6}
        \abs{\widehat{g}(\xi)}}^2
        |\xi|^{-2s}
        \int_{\frac{1}{2}|\xi| \leq |\eta| \leq 2|\xi|}
        |\xi+\eta|^{2(s-1)}
        d\eta
        d\xi}^{\frac{1}{2}}\\
        & \quad
        =
        C
        \n{f}_{\dH^{s_5}}
        \Mp{
        \int_{\mathbb{R}^2_{\xi}}
        \sp{|\xi|^{s_6}
        \abs{\widehat{g}(\xi)}}^2
        |\xi|^{-2s}
        \int_{\frac{1}{2}|\xi| \leq |\eta-\xi| \leq 2|\xi|}
        |\eta|^{2(s-1)}
        d\eta
        d\xi}^{\frac{1}{2}}\\
        & \quad
        \leq
        C
        \n{f}_{\dH^{s_5}}
        \Mp{
        \int_{\mathbb{R}^2_{\xi}}
        \sp{|\xi|^{s_6}
        \abs{\widehat{g}(\xi)}}^2
        |\xi|^{-2s}
        \int_{|\eta| \leq 3|\xi|}
        |\eta|^{2(s-1)}
        d\eta
        d\xi}^{\frac{1}{2}}\\
        & \quad
        \leq
        C
        \n{f}_{\dH^{s_5}}
        \n{g}_{\dH^{s_6}}.
    \end{align}
    For the estimate of $I_3(\xi)$ and $I_5(\xi)$, we have by $(1/2)|\eta| \leq |\xi| \leq 2|\eta|$ for $|\eta| \geq 2|\xi-\eta|$ that
    \begin{align}
        \sum_{k=3,5}
        |\xi|^{s_2-1}
        \abs{I_k(\xi)}
        \leq{}&
        \int_{|\eta| \geq 2|\xi-\eta|}
        |\eta|^{s-1}
        \abs{\widehat{f}(\eta)}
        \abs{\widehat{g}(\xi-\eta)}
        d\eta \\
        \leq{}&
        \Mp{\int_{|\eta| \geq 2|\xi-\eta|}
        \sp{
        |\eta|^{s-1}
        |\xi-\eta|^{-s_6}
        \abs{\widehat{f}(\eta)}
        }^2d\eta}^{\frac{1}{2}}
        \n{g}_{\dH^{s_6}},
    \end{align}
    which implies 
    \begin{align}
        \sum_{k=3,5}
        \n{|\xi|^{s_2-1}I_k(\xi)}_{L^2_{\xi}}^2
        \leq{}&
        C\int_{\mathbb{R}^2_{\eta}}
        \int_{|\eta| \geq 2|\xi-\eta|}
        |\xi-\eta|^{-2s_6}
        d\xi
        |\eta|^{2(s-1)}\abs{\widehat{f}(\eta)}^2
        d\eta
        \n{g}_{\dH^{s_6}}^2\\
        ={}&
        C\int_{\mathbb{R}^2_{\eta}}
        |\eta|^{-2s_6+2}
        |\eta|^{2(s-1)}\abs{\widehat{f}(\eta)}^2
        d\eta
        \n{g}_{\dH^{s_6}}^2\\
        ={}&
        C\int_{\mathbb{R}^2_{\eta}}
        |\eta|^{2s_5}\abs{\widehat{f}(\eta)}^2
        d\eta
        \n{g}_{\dH^{s_6}}^2
        =
        C
        \n{f}_{\dH^{s_5}}^2
        \n{g}_{\dH^{s_6}}^2.
    \end{align}
    For the estimate of $I_4(\xi)$, we see that 
    \begin{align}
        |\xi|^{s_2-1}
        \abs{I_4(\xi)}
        \leq{}&
        |\xi|^{s_2-1}
        \int_{\frac{1}{2}|\eta| \leq |\xi - \eta| \leq 2|\eta|}
        |\xi-\eta|^{-s_2}
        \abs{\widehat{f}(\eta)}
        |\xi-\eta|^{s}
        \abs{\widehat{g}(\xi-\eta)}
        d\eta\\
        \leq{}&
        |\xi|^{s_2-1}
        \int_{\frac{1}{2}|\eta| \leq |\xi - \eta| \leq 2|\eta|}
        |\xi-\eta|^{-s_2}
        |\eta|^{s_5}
        \abs{\widehat{f}(\eta)}
        |\xi-\eta|^{s_6}
        \abs{\widehat{g}(\xi-\eta)}
        d\eta\\
        \leq{}&
        |\xi|^{s_2-1}
        \Mp{\int_{\frac{1}{2}|\eta| \leq |\xi - \eta| \leq 2|\eta|}
        \sp{|\xi-\eta|^{-s_2}
        |\eta|^{s_5}
        \abs{\widehat{f}(\eta)}}^2
        d\eta}^{\frac{1}{2}}\n{g}_{\dH^{s_6}},
    \end{align}
    which implies 
    \begin{align}
        \n{|\xi|^{s_2-1}
        I_4(\xi)}_{L^2_{\xi}}^2
        \leq{}&
        \int_{\mathbb{R}^2}
        \sp{
        |\eta|^{s_5}
        \abs{\widehat{f}(\eta)}}^2
        \int_{\frac{1}{2}|\eta| \leq |\xi - \eta| \leq 2|\eta|}
        |\xi|^{2(s_2-1)}
        |\xi-\eta|^{-2s_2}
        d\xi
        d\eta
        \n{g}_{\dH^{s_6}}^2\\
        \leq{}&
        C
        \n{f}_{\dH^{s_5}}^2\n{g}_{\dH^{s_6}}^2,
    \end{align}
    where we have used the estimate 
    \begin{align}
        \int_{\frac{1}{2}|\eta| \leq |\xi - \eta| \leq 2|\eta|}
        |\xi|^{2(s_2-1)}
        |\xi-\eta|^{-2s_2}
        d\xi
        \leq{}&
        C|\eta|^{-2s_2}
        \int_{| \xi | \leq 3| \eta |}
        |\xi|^{2(s_2-1)}
        d\xi\\
        \leq{}&C.
    \end{align}
    Thus, we complete the proof.
\end{proof}
\begin{rem}
    In \cite{Bah-Che-Dan-11}, the first estimate of Lemma \ref{lemm:L-M} was proved via the Littlewood--Paley analysis, whereas we may provide an alternative proof by using the same arguments for the estimates of $I_2$ and $I_3$ in the above proof.
\end{rem}

\section{Proof of Theorem \ref{thm:1}}\label{sec:pf_thm:1}
In this section, we provide the proof of Theorem \ref{thm:1}.
First of all, we provide a lemma for constructing the approximation sequence of \eqref{eq:s-QG}.
Here, we introduce a notation.
For $N>0$, let $P_N$ be an operator on $L^2(\mathbb{R}^2)$ given by 
\begin{align}
    [P_N g](x)
    :={}& 
    \frac{1}{(2\pi)^2}
    \int_{|\xi| \leq 2^N} e^{ix \cdot \xi} \widehat{g}(\xi) d\xi
\end{align}
for $g \in L^2(\mathbb{R}^2)$ and $x\in \mathbb{R}^2$.
We also define
\begin{align}
    P_NL^2(\mathbb{R}^2)
    :=
    \Mp{P_Ng\ ;\ g \in L^2(\mathbb{R}^2)}.
\end{align}
\begin{lemm}\label{lemm:L-M}
    Let $0<\alpha \leq 1/2$.
    Then, there exists a positive constant $\eta_*=\eta_*(\alpha)$ such that if a given vector field $v \in P_N L^2(\mathbb{R}^2)$ satisfies $\div v= 0$ and $\n{v}_{\dH^{2-2\alpha}}\leq \eta_*$, then for every $N \in \mathbb{N}$ and $f \in \dH^{-\alpha}(\mathbb{R}^2) \cap \dH^{2-4\alpha}(\mathbb{R}^2)$, there exists a unique solution $\theta_N \in P_NL^2(\mathbb{R}^2)$ to the linear problem 
    \begin{align}\label{eq:lin}
        (-\Delta)^{\alpha} \theta_N + P_N (v \cdot \nabla \theta_N) = P_Nf.
    \end{align}
    Moreover, there exists a positive constant $C_*=C_*(\alpha)$ such that 
    \begin{align}\label{est:lin_a_priori}
        \n{\theta_N}_{\dH^{\alpha}}
        \leq 
        \n{f}_{\dH^{-\alpha}},\quad
        \n{\theta_N}_{\dH^{2-2\alpha}}
        \leq 
        C_*
        \n{f}_{\dH^{2-4\alpha}}.
    \end{align}
\end{lemm}

\begin{proof}
    Let us define a bilinear form $B_{\alpha,N}^{v}[\cdot, \cdot]$ on $P_NL^2(\mathbb{R}^2)$ as 
    \begin{align}
        B_{\alpha,N}^{v}[\theta,\varphi]
        :=
        \lr{ \theta + (-\Delta)^{-\alpha}P_N (v \cdot \nabla \theta), \varphi}_{L^2}
    \end{align}
    for $\theta, \varphi \in P_NL^2(\mathbb{R}^2)$.
    Then, from $(-\Delta)^{-\frac{\alpha}{2}}(v \cdot \nabla \theta) = (-\Delta)^{-\frac{\alpha}{2}}\div (v \theta)$, it follows that 
    \begin{align}
        \abs{B_{\alpha,N}^{v}[\theta,\varphi]}
        \leq{}&
        \n{\theta}_{L^2}
        \n{\varphi}_{L^2}
        +
        CN^{1-2\alpha}\n{v\theta}_{L^2}\n{\varphi}_{L^2}\\
        \leq{}&
        \sp{1+
        C2^{(1-2\alpha)N}
        \n{v}_{L^{\infty}}
        }
        \n{\theta}_{L^2}
        \n{\varphi}_{L^2}\\
        \leq{}&
        \sp{1+
        C2^{(4-2\alpha)N}
        \n{v}_{L^2}
        }
        \n{\theta}_{L^2}
        \n{\varphi}_{L^2}.
    \end{align}
    On the other hand,
    by the fact that the divergence free condition on $v$ yields  
    \begin{align}
        \lr{ v \cdot \nabla (-\Delta)^{-\frac{\alpha}{2}}\theta, (-\Delta)^{-\frac{\alpha}{2}}\theta}_{L^2}
        =0,
    \end{align}
    it holds that 
    \begin{align}
        \abs{\lr{(-\Delta)^{-\alpha}P_N (v \cdot \nabla \theta), \theta}_{L^2}}
        &=
        \abs{
        \lr{ \lp{v\cdot \nabla, (-\Delta)^{-\frac{\alpha}{2}}}\theta, (-\Delta)^{-\frac{\alpha}{2}}\theta }_{L^2}}\\
        &\leq
        \n{ \lp{v\cdot \nabla, (-\Delta)^{-\frac{\alpha}{2}}}\theta}_{\dH^{\alpha}}
        \n{ \theta }_{L^2}\\
        &\leq
        C_1
        \n{ v }_{\dH^{2-2\alpha}}
        \n{ \theta }_{L^2}^2
    \end{align}
    for some positive constant $C_1=C_1(\alpha)$ independent of $N$.
    Hence, assuming $\n{v}_{\dH^{2-2\alpha}} \leq 1/(2C_1)$, we have 
    \begin{align}
        \abs{B_{\alpha,N}^{v}[\theta,\theta]}
        \geq{}&
        \n{\theta}_{L^2}^2
        -
        \abs{\lr{(-\Delta)^{-\alpha}P_N (v \cdot \nabla \theta), \theta}_{L^2}}
        \\
        \geq{}&
        \n{\theta}_{L^2}^2
        -
        C_1
        \n{ v }_{\dH^{2-2\alpha}}
        \n{ \theta }_{L^2}^2\\
        \geq{}&
        \frac{1}{2}
        \n{\theta}_{L^2}^2.
    \end{align}
    Hence, by the Lax–Milgram theorem, there exists a unique $\theta_N \in P_NL^2(\mathbb{R}^2)$ such that $B_{\alpha,N}^{v}[\theta,\varphi] = \lr{(-\Delta)^{-\alpha}P_Nf,\varphi}_{L^2}$ for every $\varphi \in P_NL^2(\mathbb{R}^2)$, which completes the proof of the unique existence of the solution.

    Next, we prove the a priori estimate \eqref{est:lin_a_priori}.
    For the $\dH^{\alpha}(\mathbb{R}^2)$-estimate, taking the $L^2(\mathbb{R}^2)$-inner product of the equation with $\theta_N$, we have 
    \begin{align}
        \n{(-\Delta)^{\frac{\alpha}{2}}\theta_N}_{L^2}^2
        =
        \lr{f,\theta_N}_{L^2}
        \leq
        \n{f}_{\dH^{-\alpha}}
        \n{\theta_N}_{\dH^{\alpha}},
    \end{align}
    where we have used $\lr{P_N(v \cdot \nabla \theta_N),\theta_N}_{L^2}=\lr{v \cdot \nabla \theta_N,\theta_N}_{L^2}=0$ due to $\div v=0$.
    Thus, we have
    \begin{align}
        \n{\theta_N}_{\dH^{\alpha}}
        \leq
        \n{f}_{\dH^{-\alpha}}.
    \end{align}
    For the $\dH^{2-2\alpha}(\mathbb{R}^2)$-estimate, taking the inner $\dH^{2-3\alpha}(\mathbb{R}^2)$-inner product with $\theta_N$, we have 
    \begin{align}
        \n{\theta_N}_{\dH^{2-2\alpha}}^2
        +{}
        \lr{v\cdot \nabla \theta_N, \theta_N}_{\dH^{2-3\alpha}}
        ={}
        \lr{f,\theta_N}_{\dH^{2-3\alpha}},
    \end{align}
    It follows from $\div v = 0$ that
    \begin{align}
        \lr{v\cdot \nabla \theta_N,\theta_N}_{\dH^{2-3\alpha}}
        =
        -
        \lr{\lp{v\cdot\nabla, (-\Delta)^{\frac{2-3\alpha}{2}}}\theta_N,(-\Delta)^{\frac{2-3\alpha}{2}}\theta_N}_{L^2}.
    \end{align}
    Moreover, we see that
    \begin{align}
        \abs{\lr{v\cdot \nabla \theta_N,\theta_N}_{\dH^{2-3\alpha}}
        }
        &
        ={}
        \abs{\lr{\lp{v\cdot \nabla, (-\Delta)^{\frac{2-3\alpha}{2}}} \theta_N,(-\Delta)^{\frac{2-3\alpha}{2}}\theta_N}_{L^2}
        }
        \\
        &
        \leq{}
        \abs{\lr{\lp{v\cdot \nabla, (-\Delta)^{\frac{2-3\alpha}{2}}} \theta_N,(-\Delta)^{\frac{2-3\alpha}{2}}\theta_N}_{L^2}
        }
        \\
        &
        \leq{}
        \n{\lp{v\cdot \nabla, (-\Delta)^{\frac{2-3\alpha}{2}}} \theta_N}_{\dH^{-\alpha}}
        \n{\theta_N}_{\dH^{2-2\alpha}}
        \\
        &
        \leq{}
        C
        \n{v}_{\dH^{2-2\alpha}}
        \n{\theta_N}_{\dH^{2-2\alpha}}^2
    \end{align}
    and 
    $\abs{
    \lr{f,\theta_N}_{\dH^{2-3\alpha}}
    }
    \leq{} 
    \n{f}_{\dH^{2-4\alpha}}
    \n{\theta_N}_{\dH^{2-2\alpha}}$.
    Thus, we have
    \begin{align}
        \n{\theta_N}_{\dH^{2-2\alpha}}
        \leq{}&
        C_2\n{f}_{\dH^{2-4\alpha}}
        +
        C_2\n{v}_{\dH^{2-2\alpha}}\n{\theta_N}_{\dH^{2-2\alpha}}
    \end{align}
    for some positive constant $C_2=C_2(\alpha)$.
    Hence, assuming that $v$ satisfy 
    \begin{align}
        \n{v}_{\dH^{2-2\alpha}}
        \leq
        \min \Mp{\frac{1}{2C_1},\frac{1}{2C_2}}=:\eta_*,
    \end{align}
    we complete the proof.
\end{proof}

Now we are in a position to present the proof of Theorem \ref{thm:1}.

\begin{proof}[Proof of Theorem \ref{thm:1}]
    Let $C_*$ and $\eta_*$ be the positive constant to be determined in Lemma \ref{lemm:L-M} and assume that $f \in \dH^{\alpha}(\mathbb{R}^2) \cap \dH^{2-2\alpha}(\mathbb{R}^2)$ satisfies
    \begin{align}
        \n{f}_{\dH^{2-2\alpha}} \leq \delta_{**}:=\frac{\eta_{**}}{C_*},
    \end{align}
    where the constant $0 < \eta_{**} \leq \eta_*$ is to be determined later.
    Let us define the approximation sequence $\{\theta_N\}_{N=0}^{\infty}$ for \eqref{eq:s-QG} inductively as follows:
    Let $\theta_0:=0$ and $\theta_1:=(-\Delta)^{-\alpha}P_1f$.
    Next, we suppose that an appropriate $\theta_N \in P_N L^2(\mathbb{R}^2)$ with $\n{\theta_N}_{\dH^{2-2\alpha}} \leq \delta_*$ is determined and let $v_N:=\nabla^{\perp}(-\Delta)^{-1/2}\theta_N$.
    Then, using Lemma \ref{lemm:L-M} via the estimate $\n{v_N}_{\dH^{2-2\alpha}} \leq \n{\theta_N}_{\dH^{2-2\alpha}}\leq \delta_*$, we may define $\theta_{N+1} \in P_{N+1}L^2(\mathbb{R}^2)$ as a unique solution to the equation
    \begin{align}\label{eq:app_s-QG}
        (-\Delta)^{\alpha} \theta_{N+1} + P_{N+1}( v_N \cdot \nabla \theta_{N+1}) = P_{N+1}f
    \end{align}
    with the estimate 
    \begin{align}
        \n{\theta_{N+1}}_{\dH^{\alpha}}
        \leq 
        \n{f}_{\dH^{-\alpha}},\quad
        \n{\theta_{N+1}}_{\dH^{2-2\alpha}}
        \leq 
        C_*
        \n{f}_{\dH^{2-4\alpha}}
        \leq 
        \eta_{**}.
    \end{align}
    Then, the sequence $\{ \theta_N \}_{N=0}^{\infty}$ has the uniform estimate
    \begin{align}\label{est:unif}
        \sup_{N \in \mathbb{N}}
        \n{\theta_{N}}_{\dH^{\alpha}}
        \leq 
        \n{f}_{\dH^{-\alpha}},\quad
        \sup_{N \in \mathbb{N}}
        \n{\theta_{N}}_{\dH^{2-2\alpha}}
        \leq
        C_*
        \n{f}_{\dH^{2-4\alpha}}\leq
        \eta_{**}.
    \end{align}
    To deduce the convergence of $\{ \theta_N \}_{N=1}^{\infty}$, we consider the difference $\theta_{N+1} - \theta_N$.
    Since, it holds
    \begin{align}
        (-\Delta)^{\alpha}(\theta_{N+1}-\theta_N)
        &
        +
        P_{N+1}
        (v_N \cdot \nabla (\theta_{N+1}-\theta_N))\\
        &
        +
        P_{N+1}
        ((v_N -v_{N-1})\cdot \nabla \theta_N)\\
        &
        +
        (P_{N+1}-P_N)
        (v_{N-1} \cdot \nabla \theta_N)
        =(P_{N+1}-P_N)f.
    \end{align}
    Taking $L^2(\mathbb{R}^2)$-inner product of the above equation with $\theta_{N+1}-\theta_N$, we have by similar calculation as above that
    \begin{align}
        \n{\theta_{N+1}-\theta_N}_{\dH^{\alpha}}
        \leq{}&
        \n{(P_{N+1}-P_N)f}_{\dH^{-\alpha}}\\
        &
        +
        \n{(P_{N+1}-P_N)
        (v_{N-1}\cdot \nabla \theta_N)}_{\dH^{-\alpha}}
        \\
        &
        +
        \n{P_{N+1}
        ((v_N -v_{N-1})\cdot \nabla \theta_N)}_{\dH^{-\alpha}}.
    \end{align}
    Since all terms in the right hand side are bounded as 
    \begin{align}
        &
        \begin{aligned}
        \n{(P_{N+1}-P_N)f}_{\dH^{-\alpha}}
        \leq{}
        C
        2^{-(2-3\alpha)N}
        \n{f}_{\dH^{2-4\alpha}}
        \leq{}
        C
        2^{-\frac{\alpha}{2}N}
        \n{f}_{\dH^{2-4\alpha}},
        \end{aligned}
        \\
        &
        \begin{aligned}
        \n{(P_{N+1}-P_N)
        (v_{N-1}\cdot \nabla \theta_N)}_{\dH^{-\alpha}}
        \leq{}&
        C2^{-\frac{\alpha}{2} N}
        \n{v_{N-1} \cdot \nabla \theta_N}_{\dH^{-\frac{\alpha}{2}}}\\
        \leq{}&
        C2^{-\frac{\alpha}{2} N}
        \n{v_{N-1}}_{\dH^{\frac{3}{2}\alpha}}
        \n{\theta_N}_{\dH^{2-2\alpha}}\\
        \leq{}&
        C2^{-\frac{\alpha}{2} N}
        \n{f}_{\dH^{-\alpha} \cap \dH^{2-4\alpha}}^2,
        \end{aligned}\\
        &
        \begin{aligned}
        \n{P_{N+1}
        ((v_N-v_{N-1})\cdot \nabla \theta_N)}_{\dH^{-\alpha}}
        \leq{}&
        C
        \n{\theta_N}_{\dH^{2-2\alpha}}
        \n{\theta_{N}-\theta_{N-1}}_{\dH^{\alpha}}\\
        \leq{}&
        C
        \eta_{**}
        \n{\theta_{N}-\theta_{N-1}}_{\dH^{\alpha}},
        \end{aligned}
    \end{align}
    we see that 
    \begin{align}
        \n{\theta_{N+1}-\theta_N}_{\dH^{\alpha}}
        \leq{}&
        C_{**}
        2^{- \frac{\alpha}{2} N}
        \sp{
        \n{f}_{\dH^{2-4\alpha}}
        +
        \n{f}_{\dH^{-\alpha} \cap \dH^{2-4\alpha}}^2
        }\\
        &
        +
        C_{**}
        \eta_{**}
        \n{\theta_{N}-\theta_{N-1}}_{\dH^{\alpha}}.
    \end{align}
    with some positive constant $C_{**}=C_{**}(\alpha)$.
    Then, choosing $\eta_{**}:=\min\{\eta_*,1/(2C_{**})\}$, we see that 
    \begin{align}
        \n{\theta_{N+1}-\theta_N}_{\dH^{\alpha}}
        \leq{}&
        C_{**}
        2^{- \frac{\alpha}{2} N}
        \sp{
        \n{f}_{\dH^{-\alpha} \cap \dH^{2-4\alpha}}
        +
        \n{f}_{\dH^{-\alpha} \cap \dH^{2-4\alpha}}^2
        }\\
        &
        +
        \frac{1}{2}
        \n{\theta_{N}-\theta_{N-1}}_{\dH^{\alpha}},
    \end{align}
    which implies $\{ \theta_N \}_{N=1}^{\infty}$ is a Cauchy sequence in $\dH^{\alpha}(\mathbb{R}^2)$.
    Hence, there exists a limit $\theta_N \xrightarrow{N \to \infty} \theta$ strongly in $\dH^{\alpha}(\mathbb{R}^2)$ and the uniform estimate \eqref{est:unif} ensures $\theta \in \dH^{\alpha} (\mathbb{R}^2) \cap \dH^{2-2\alpha}(\mathbb{R}^2)$ with the estimate 
    \begin{align}
        \n{\theta}_{\dH^{\alpha}}
        \leq 
        \n{f}_{\dH^{-\alpha}},\quad
        \n{\theta}_{\dH^{2-2\alpha}}
        \leq 
        C_*
        \n{f}_{\dH^{2-4\alpha}}\leq \eta_{**}.
    \end{align}
    We may easily see that this $\theta$ solves \eqref{eq:s-QG} in the distributional sense.
    We omit the proof of the uniqueness of the solution since it is similar to the calculation performed in the convergence part of the approximation sequence.
\end{proof}

\section{Proof of Theorem \ref{thm:2}}\label{sec:pf_thm:2}
In this section, we prove Theorem \ref{thm:2}.
We begin with the following lemma for the smoothing operator $e^{\varepsilon^2\Delta}:=\mathscr{F}^{-1}e^{-\varepsilon^2 |\xi|^2}\mathscr{F}$.
\begin{lemm}\label{lemm:lim}
    Let $s \in \mathbb{R}$ and $0\leq \sigma<2$.
    For every $f \in \dH^{s}(\mathbb{R}^2)$, it holds
    \begin{align}\label{intpln}
        \varepsilon^{-\sigma}
        \n{e^{\varepsilon^2\Delta}f-f}_{\dH^{s-\sigma}}
        \leq 
        \sqrt{2}\n{f}_{\dH^s}^{\frac{\sigma}{2}}\n{e^{\varepsilon^2\Delta}f-f}_{\dH^s}^{1-\frac{\sigma}{2}}.
    \end{align}
    In particular, it holds 
    \begin{align}\label{cor-lim}
        \varepsilon^{-\sigma}
        \n{e^{\varepsilon^2\Delta}f-f}_{\dH^{s-\sigma}}
        \leq 
        2^{\frac{3-\sigma}{2}}\n{f}_{\dH^s}, \quad
        \lim_{\varepsilon \downarrow 0}
        \varepsilon^{-\sigma}
        \n{e^{\varepsilon^2\Delta}f-f}_{\dH^{s-\sigma}}
        =0.
    \end{align}
    Moreover, for every compact set $K \subset \dH^{s}(\mathbb{R}^2)$, it holds 
    \begin{align}\label{unif-lim}
        \lim_{\varepsilon \downarrow 0}
        \sup_{f \in K}
        \varepsilon^{-\sigma}
        \n{e^{\varepsilon^2 \Delta} f-f}_{\dH^{s-\sigma}}
        =
        0.
    \end{align}
\end{lemm}
\begin{proof}
    Let us first show \eqref{intpln}.
    It suffices to consider the case $f \neq 0$.
    For $R>0$, we see that 
    \begin{align}
        \varepsilon^{-2\sigma}
        \n{e^{\varepsilon^2\Delta}f-f}_{\dH^{s-\sigma}}^2
        ={}&
        \int_{\mathbb{R}^2}
        \sp{\frac{1-e^{\varepsilon^2|\xi|^2}}{(\varepsilon|\xi|)^{\sigma}}}^2
        |\xi|^{2s}\abs{\widehat{f}(\xi)}^2
        d\xi\\
        \leq{}&
        \int_{|\xi|\leq R}
        (\varepsilon|\xi|)^{2(2-\sigma)}
        |\xi|^{2s}\abs{\widehat{f}(\xi)}^2
        d\xi\\
        &
        +
        \int_{|\xi|\geq R}
        (\varepsilon|\xi|)^{-2\sigma}
        \sp{1-e^{-\varepsilon^2|\xi|^2}}
        |\xi|^{2s}\abs{\widehat{f}(\xi)}^2
        d\xi\\
        \leq{}&
        (\varepsilon R)^{2(2-\sigma)}
        \n{f}_{\dH^s}^2
        +
        (\varepsilon R)^{-2\sigma}
        \n{e^{\varepsilon^2\Delta}f-f}_{\dH^{s}}^2.
    \end{align}
    Then, to obtain \eqref{intpln}, it suffices to take $R$ so that 
    \begin{align}
        (\varepsilon R)^{2-\sigma}
        \n{f}_{\dH^s}
        =
        (\varepsilon R)^{-\sigma}
        \n{e^{\varepsilon^2\Delta}f-f}_{\dH^{s}}.
    \end{align}
    Using \eqref{intpln}, the boundedness of the heat kernel on homogeneous Sobolev spaces, and the dominated convergence theorem, we immediately obtain \eqref{cor-lim}.
    
    Next, we prove the uniformity of the limit \eqref{unif-lim}.
    Let $\delta>0$ be arbitrary.
    Then, by the compactness, there exists a finitely many $f_1,\dots, f_M \in K$ such that 
    \begin{align}
        K \subset \bigcup_{m=1}^M B(f_m,\delta),
    \end{align}
    where $B(f,\delta)$ be the open ball in $\dH^s(\mathbb{R}^2)$ centered at $f$ with the radius $\delta$.
    For any $f \in K$, we may choose a $m_f \in \{1,\dots,M\}$ so that $f \in B(f_{m_f},\delta)$.
    Then, we see that 
    \begin{align}
        \varepsilon^{-\sigma}
        \n{e^{\varepsilon^2 \Delta} f-f}_{\dH^{s-\sigma}}
        \leq {}&
        \varepsilon^{-\sigma}
        \n{e^{\varepsilon^2 \Delta} (f-f_{m_f}) -(f-f_{m_f})}_{\dH^{s-\sigma}}\\
        &
        + 
        \varepsilon^{-\sigma}
        \n{e^{\varepsilon^2 \Delta} f_{m_f}-f_{m_f}}_{\dH^{s-\sigma}}\\
        \leq {}&
        2^{\frac{3-\sigma}{2}}\n{f-f_{m_f}}_{\dH^s} 
        + 
        \varepsilon^{-\sigma}\n{e^{\varepsilon^2 \Delta} f_{m_f}-f_{m_f}}_{\dH^{s-\sigma}}\\
        \leq {}&
        2^{\frac{3-\sigma}{2}}\delta + \max _{m=1,\dots,M}
        \varepsilon^{-\sigma}
        \n{e^{\varepsilon^2 \Delta} f_{m}-f_{m}}_{\dH^{s-\sigma}}.
    \end{align}
    This implies 
    \begin{align}
        \limsup_{\varepsilon \downarrow 0}
        \sup_{f \in K}
        \varepsilon^{-\sigma}
        \n{e^{\varepsilon^2 \Delta} f-f}_{\dH^{s-\sigma}}
        \leq {}
        2^{\frac{3-\sigma}{2}}
        \delta.
    \end{align}
    Letting $\delta \downarrow 0$, we complete the proof
\end{proof}
Now, we are ready to show Theorem \ref{thm:2}.
\begin{proof}[Proof of Theorem \ref{thm:2}]
The continuous dependence in $\dH^{\alpha}(\mathbb{R}^2)$ is proved similarly to the proof of the convergence of the approximation sequence in Section \ref{sec:pf_thm:1}.
Thus, we focus on the continuity in $\dH^{2-2\alpha}(\mathbb{R}^2)$.
Our proof is based on the approximation of smooth solutions that is inspired by classical results Bona--Smith \cite{BS} and Kato--Lai \cite{KL}.
In the following, let $n \in \bar{\mathbb{N}}$.
Let us consider the solution $\theta_n^{\varepsilon}$ to the stationary quasi-geostrophic equation with the regularized external force $e^{\varepsilon^2\Delta}f_n$:
\begin{align}\label{eq:app_ep}
    \begin{cases}
        (-\Delta)^{\alpha} \theta_n^{\varepsilon} + v_n^{\varepsilon} \cdot \nabla \theta_n^{\varepsilon} = e^{\varepsilon^2\Delta}f_n, \\
        v_n^{\varepsilon} = \nabla^{\perp}(-\Delta)^{-1/2}\theta_n^{\varepsilon}.
    \end{cases}
\end{align}
Then, from the same energy calculation in Section \ref{sec:pf_thm:1}, we see that 
\begin{align}\label{est:unif_n}
    \sup_{\varepsilon > 0}
    \n{\theta_n^{\varepsilon}}_{\dH^{\alpha}}
    \leq 
    \n{f_n}_{\dH^{-\alpha}},\quad
    \sup_{\varepsilon > 0}
    \n{\theta_n^{\varepsilon}}_{\dH^{2-2\alpha}}
    \leq
    C_*
    \n{f_n}_{\dH^{2-4\alpha}}\leq
    \eta_{**}
\end{align}
for all $n \in \bar{\mathbb{N}}$.
Since $\theta_n-\theta_n^{\varepsilon}$ solves
\begin{align}\label{eq:dist}
    (-\Delta)^{\alpha}(\theta_n-\theta_n^{\varepsilon})
    +
    v_n \cdot \nabla (\theta_n-\theta_n^{\varepsilon})
    +
    (v_n - v_n^{\varepsilon}) \cdot \nabla \theta_n^{\varepsilon}
    =
    f_n
    -
    e^{\varepsilon^2\Delta}f_n,
\end{align}
we see by taking the $\dH^{2-3\alpha}(\mathbb{R}^2)$-inner product of the second equation of \eqref{eq:dist} with $\theta_n-\theta_n^{\varepsilon}$ and calculating the similar energy argument as in the previous section that
\begin{align}
    \n{\theta_n-\theta_n^{\varepsilon}}_{\dH^{2-2\alpha}}
    \leq{}&
    \n{f_n - e^{\varepsilon^2\Delta}f_n}_{\dH^{2-4\alpha}}\\
    &
    +
    \n{\lp{(-\Delta)^{\frac{2-3\alpha}{2}}, v_n \cdot \nabla}(\theta_n^{\varepsilon}-\theta_n)}_{\dH^{-\alpha}}
    +
    \n{(v_n^{\varepsilon} - v_n) \cdot \nabla \theta_n^{\varepsilon}}_{\dH^{2-4\alpha}}\\
    \leq{}&
    \n{f_n - e^{\varepsilon^2\Delta}f_n}_{\dH^{2-4\alpha}}\\
    &
    +
    C
    \n{\theta_n}_{\dH^{2-2\alpha}}\n{\theta_n-\theta_n^{\varepsilon}}_{\dH^{2-2\alpha}}
    +
    C
    \n{\theta_n^{\varepsilon} - \theta_n}_{\dH^{\alpha}}\n{\theta_n^{\varepsilon}}_{\dH^{4-5\alpha}}.
\end{align}
For the estimate of $\n{\theta_n^{\varepsilon}}_{\dH^{4-5\alpha}}$, taking $\dH^{4-6\alpha}(\mathbb{R}^2)$-inner product of \eqref{eq:app_ep} with $\theta_n^{\varepsilon}$, we have 
\begin{align}
    \n{\theta_n^{\varepsilon}}_{\dH^{4-5\alpha}}^2
    \leq {}&
    \n{e^{\varepsilon^2\Delta}f_n}_{\dH^{4-7\alpha}}
    \n{\theta_n^{\varepsilon}}_{\dH^{4-5\alpha}}
    +
    \abs{
    \lr{v_n^{\varepsilon}\cdot \nabla\theta_n^{\varepsilon},\theta_n^{\varepsilon}}_{\dH^{4-6\alpha}}}\\
    = {}&
    \n{e^{\varepsilon^2\Delta}f_n}_{\dH^{4-7\alpha}}
    \n{\theta_n^{\varepsilon}}_{\dH^{4-5\alpha}}
    +
    \abs{
    \lr{\lp{(-\Delta)^{\frac{4-6\alpha}{2}}, v_n^{\varepsilon}\cdot \nabla}\theta_n^{\varepsilon},\theta_n^{\varepsilon}}_{\dH^{4-6\alpha}}}.
\end{align}
Thus, we see that 
\begin{align}
    \n{\theta_n^{\varepsilon}}_{\dH^{4-5\alpha}}
    \leq{}&
    \n{e^{\varepsilon^2\Delta}f_n}_{\dH^{4-7\alpha}}
    +
    \n{\lp{(-\Delta)^{\frac{4-6\alpha}{2}}, v_n^{\varepsilon}\cdot \nabla}\theta_n^{\varepsilon}}_{\dH^{-\alpha}}
    \\
    \leq{}&
    \n{e^{\varepsilon^2\Delta}f_n}_{\dH^{4-7\alpha}}
    +
    C
    \n{v_n^{\varepsilon}}_{\dH^{2-2\alpha}}
    \n{\theta_n^{\varepsilon}}_{\dH^{4-5\alpha}}
    +
    C
    \n{v_n^{\varepsilon}}_{\dH^{4-5\alpha}}
    \n{\theta_n^{\varepsilon}}_{\dH^{2-2\alpha}}\\
    \leq{}&
    \n{e^{\varepsilon^2\Delta}f_n}_{\dH^{4-7\alpha}}
    +
    C
    \n{\theta_n^{\varepsilon}}_{\dH^{2-2\alpha}}
    \n{\theta_n^{\varepsilon}}_{\dH^{4-5\alpha}},
\end{align}
which and \eqref{est:unif_n} imply 
\begin{align}
    \n{\theta_n^{\varepsilon}}_{\dH^{4-5\alpha}}
    \leq{}&
    C\n{e^{\varepsilon^2\Delta}f_n}_{\dH^{4-7\alpha}}.
\end{align}
For the estimate of $\n{\theta_n-\theta_n^{\varepsilon}}_{\dH^{\alpha}}$, taking the $L^2(\mathbb{R}^2)$-inner product of the second equation of \eqref{eq:dist} with $\theta_n^{\varepsilon}-\theta_n$, we have
\begin{align}
    \n{\theta_n-\theta_n^{\varepsilon}}_{\dH^{\alpha}}
    \leq{}&
    \n{f_n - e^{\varepsilon^2\Delta} f_n}_{\dH^{-\alpha}}
    +
    \n{(v_n - v_n^{\varepsilon}) \cdot \nabla \theta_n^{\varepsilon}}_{\dH^{-\alpha}}\\
    \leq{}&
    \n{f_n - e^{\varepsilon^2\Delta}f_n}_{\dH^{-\alpha}}
    +
    C
    \n{\theta_n^{\varepsilon} - \theta_n}_{\dH^{\alpha}}
    \n{\theta_n^{\varepsilon}}_{\dH^{2-2\alpha}},
\end{align}
which and \eqref{est:unif_n} imply 
\begin{align}\label{est:dist_low}
    \n{\theta_n-\theta_n^{\varepsilon}}_{\dH^{\alpha}}
    \leq{}&
    C
    \n{f_n - e^{\varepsilon^2\Delta} f_n}_{\dH^{-\alpha}}.
\end{align}
Hence, we obtain
\begin{align}
    \n{\theta_n^{\varepsilon}-\theta_n}_{\dH^{2-2\alpha}}
    \leq{}&
    \n{f_n - e^{\varepsilon^2\Delta}f_n}_{\dH^{2-4\alpha}}\\
    &
    +
    C
    \n{e^{\varepsilon^2\Delta}f_n}_{\dH^{4-7\alpha}}
    \n{f_n - e^{\varepsilon^2\Delta} f_n}_{\dH^{-\alpha}}\\
    &
    +
    C
    \n{\theta_n}_{\dH^{2-2\alpha}}\n{\theta_n^{\varepsilon}-\theta_n^{\varepsilon}}_{\dH^{2-2\alpha}}.
\end{align}
Here, it follows from \eqref{est:unif_n} and 
$\n{e^{\varepsilon^2\Delta}f_n}_{\dH^{4-7\alpha}}
\leq{}
C
\varepsilon^{-(2-3\alpha)}
\n{f_n}_{\dH^{2-4\alpha}}$
that
\begin{align}
    \n{\theta_n^{\varepsilon}-\theta_n}_{\dH^{2-2\alpha}}
    \leq{}&
    \n{f_n - e^{\varepsilon^2\Delta}f_n}_{\dH^{2-4\alpha}}
    +
    C\varepsilon^{-(2-3\alpha)}
    \n{f_n - e^{\varepsilon^2\Delta} f_n}_{\dH^{-\alpha}}
    \n{f_n}_{\dH^{-\alpha}}.
\end{align}
On the other hand, since it holds
\begin{align}\label{est:dist-2}
    (-\Delta)^{\alpha}(\theta_n^{\varepsilon}-\theta_\infty^{\varepsilon})
    +
    v_n \cdot \nabla (\theta_n^{\varepsilon}-\theta_\infty^{\varepsilon})
    +
    (v_n^{\varepsilon} - v_\infty^{\varepsilon}) \cdot \nabla \theta_\infty^{\varepsilon}
    =e^{\varepsilon^2\Delta}(f_n-f_{\infty}),
\end{align}
the energy argument yields 
\begin{align}
    \n{\theta_n^{\varepsilon}-\theta_\infty^{\varepsilon}}_{\dH^{\alpha}}
    &
    \leq
    \n{e^{\varepsilon^2\Delta}(f_n-f_{\infty})}_{\dH^{-\alpha}}
    +
    \n{(v_n^{\varepsilon} - v_\infty^{\varepsilon}) \cdot \nabla \theta_\infty^{\varepsilon}}_{\dH^{-\alpha}}
    \\
    &
    \leq
    \n{f_n-f_{\infty}}_{\dH^{-\alpha}}
    +
    C
    \n{\theta_n^{\varepsilon} -\theta_\infty^{\varepsilon}}_{\dH^{\alpha}}
    \n{\theta_\infty^{\varepsilon}}_{\dH^{2-2\alpha}}.
\end{align}
Hence, by \eqref{est:unif_n}, we obtain 
\begin{align}
    \n{\theta_n^{\varepsilon}-\theta_\infty^{\varepsilon}}_{\dH^{\alpha}}
    \leq
    C
    \n{f_n-f_{\infty}}_{\dH^{-\alpha}}.
\end{align}
Similarly, we see by the energy calculation that
\begin{align}
    \n{\theta_n^{\varepsilon}-\theta_{\infty}^{\varepsilon}}_{\dH^{2-2\alpha}}
    \leq{}&
    \n{e^{\varepsilon^2\Delta}(f_n-f_{\infty})}_{\dH^{2-4\alpha}}\\
    &
    +
    \n{\lp{(-\Delta)^{\frac{2-3\alpha}{2}}, v_n \cdot \nabla}(\theta_n^{\varepsilon}-\theta_{\infty}^{\varepsilon})}_{\dH^{-\alpha}}\\
    &
    +
    \n{(v_n^{\varepsilon} - v_{\infty}^{\varepsilon}) \cdot \nabla \theta_n^{\varepsilon}}_{\dH^{2-4\alpha}}\\
    \leq{}&
    \n{f_n-f_{\infty}}_{\dH^{2-4\alpha}}\\
    &
    +
    C
    \n{\theta_n^{\varepsilon}}_{\dH^{2-2\alpha}}\n{\theta_n^{\varepsilon}-\theta_{\infty}^{\varepsilon}}_{\dH^{2-2\alpha}}
    \\
    &
    +
    C
    \n{\theta_n^{\varepsilon} - \theta_{\infty}^{\varepsilon}}_{\dH^{\alpha}}\n{\theta_{\infty}^{\varepsilon}}_{\dH^{4-5\alpha}}
    +
    C
    \n{\theta_n^{\varepsilon} - \theta_{\infty}^{\varepsilon}}_{\dH^{2-2\alpha}}\n{\theta_{\infty}^{\varepsilon}}_{\dH^{2-2\alpha}}.
\end{align}
Using the smallness for $\n{\theta_n}_{\dH^{2-2\alpha}}$ and $\n{\theta_n^{\varepsilon}}_{\dH^{2-2\alpha}}$ and
\begin{align}
    \n{\theta_n^{\varepsilon} - \theta_{\infty}^{\varepsilon}}_{\dH^{\alpha}}
    \n{\theta_{\infty}^{\varepsilon}}_{\dH^{4-5\alpha}}
    \leq{}&
    C
    \n{f_n-f_{\infty}}_{\dH^{-\alpha}}
    \n{e^{\varepsilon^2\Delta}f_{\infty}}_{\dH^{4-7\alpha}},
\end{align}
we obtain 
\begin{align}
    \n{\theta_n^{\varepsilon}-\theta_{\infty}^{\varepsilon}}_{\dH^{2-2\alpha}}
    \leq{}&
    C
    \n{f_n-f_{\infty}}_{\dH^{2-4\alpha}}\\
    &
    +
    C
    \n{f_n-f_{\infty}}_{\dH^{-\alpha}}
    \n{e^{\varepsilon^2\Delta}f_{\infty}}_{\dH^{4-7\alpha}},
\end{align}
which implies that for every $\varepsilon>0$, it holds
\begin{align}
    \lim_{n \to \infty}
    \n{\theta_n^{\varepsilon}-\theta_{\infty}^{\varepsilon}}_{\dH^{2-2\alpha}}
    =0.
\end{align}
Hence, there holds
\begin{align}
    \n{\theta_n- \theta_{\infty}}_{\dH^{2-2\alpha}}
    \leq{}&
    \n{\theta_n- \theta_n^{\varepsilon}}_{\dH^{2-2\alpha}}
    +
    \n{\theta_n^{\varepsilon} - \theta_{\infty}^{\varepsilon}}_{\dH^{2-2\alpha}}
    +
    \n{\theta_\infty^{\varepsilon} - \theta_{\infty}}_{\dH^{2-2\alpha}}\\
    \leq{}&
    2
    \sup_{m \in \bar{\mathbb{N}}}
    \n{\theta_m- \theta_m^{\varepsilon}}_{\dH^{2-2\alpha}}
    +
    \n{\theta_n^{\varepsilon} - \theta_{\infty}^{\varepsilon}}_{\dH^{2-2\alpha}}
    \\
    \leq{}&
    C
    \sup_{m \in \bar{\mathbb{N}}}
    \n{f_m - e^{\varepsilon^2\Delta}f_m}_{\dH^{2-4\alpha}}
    \\
    &
    +
    C
    \sup_{m \in \bar{\mathbb{N}}}
    \varepsilon^{-(2-3\alpha)}
    \n{f_m - e^{\varepsilon^2\Delta} f_m}_{\dH^{-\alpha}}\\
    &
    +
    C
    \n{\theta_n^{\varepsilon}-\theta_{\infty}^{\varepsilon}}_{\dH^{2-2\alpha}}.
\end{align}
This and Lemma \ref{lemm:lim} yield
\begin{align}
    \limsup_{n \to \infty}
    \n{\theta_n- \theta_{\infty}}_{\dH^{2-2\alpha}}
    \leq{}&
    C
    \sup_{m \in \bar{\mathbb{N}}}
    \n{f_m - e^{\varepsilon^2\Delta}f_m}_{\dH^{2-4\alpha}}
    \\
    &
    +
    C
    \sup_{m \in \bar{\mathbb{N}}}
    \varepsilon^{-(2-3\alpha)}
    \n{f_m - e^{\varepsilon^2\Delta} f_m}_{\dH^{-\alpha}}
    \to 0 \quad {\rm as} \quad \varepsilon \downarrow 0
\end{align}
and we complete the proof.
\end{proof}

\section{Proof of Theorem \ref{thm:3}}\label{sec:pf_thm:3}
The aim of this section is to prove Theorem \ref{thm:3}.
\begin{proof}[Proof of Theorem \ref{thm:3}]
We separate the proof into three steps.
\subsection*{Step 1. Definition of the sequences of external forces}
Let $\phi \in \mathscr{S}(\mathbb{R})$ so that $\widehat{\phi}$ is a even function satisfying $0 < \widehat{\phi}(\tau)\leq 1$, $\widehat{\phi}(\tau) = 1$ for $|\tau| \leq 1$, and $\widehat{\phi}(\tau) = 0$ for $|\tau| \geq 2$.
We define $f_n$ and $g_n$ as
\begin{align}
    f_n(x) &: = g_n(x) + h_n(x),\\
    g_n(x) &: = \delta2^{-(2-2\alpha)n}(-\Delta)^{\alpha}\lp{\phi(x_1)\phi(x_2)\sin (2^nx_1)},\\
    h_n(x) &: = \delta2^{-(1-2\alpha)n}(-\Delta)^{\alpha+\frac{1}{2}}\lp{\phi(x_1)\phi(x_2)}, 
\end{align}
where $\delta$ is a positive constant to be determined later.
Then, we see that  
\begin{align}
    \n{g_n}_{\dH^{-\alpha}\cap\dH^{2-4\alpha}} \leq C\delta
\end{align}
and
\begin{align}
    \n{f_n-g_n}_{\dH^{-\alpha}\cap\dH^{2-4\alpha}}
    =
    \n{h_n}_{\dH^{\alpha}\cap\dH^{2-2\alpha}}
    \leq
    C\delta2^{-(1-2\alpha)n} \to 0
\end{align}
as $n \to \infty$.
Choosing $\delta$ sufficiently small, then $f_n$ and $g_n$ generate the solutions $\theta[f_n]$ and $\theta[g_n]$ to \eqref{eq:s-QG}, respectively.
Note that these solutions satisfy 
\begin{align}
    \n{\theta[f_n]}_{\dH^{\alpha}\cap\dH^{2-2\alpha}}\leqslant C\delta,\quad
    \n{\theta[g_n]}_{\dH^{\alpha}\cap\dH^{2-2\alpha}}\leqslant C\delta.
\end{align}
\subsection*{Step 2. Estimates for the first and second iterations}
Let us define the first and second iteration operators by
\begin{align}
    \theta_1[a]:=(-\Delta)^{-\alpha}a, 
    \quad
    \theta_2[a]:=-\mathcal{B}[a,a],
    \quad
    a \in \mathscr{S}(\mathbb{R}^2)
\end{align}
where we have set 
\begin{align}
    \mathcal{B}[a,b]
    :=
    (-\Delta)^{-\alpha}[\nabla^{\perp}(-\Delta)^{-1/2}\theta_1[a]\cdot \nabla \theta_1[b]] ,\quad
    a,b \in \mathscr{S}(\mathbb{R}^2).
\end{align}
Then, we immediately see that 
\begin{align}
    \n{\theta_1[f_n]-\theta_1[g_n]}_{\dH^{2-2\alpha}}
    \leq
    C
    \n{f_n-g_n}_{\dH^{-\alpha}\cap\dH^{2-4\alpha}}
    \to 0 \quad {\rm as}\quad n \to \infty.
\end{align}
Next, we focus on the estimates of the second iterations.
We then see that
\begin{align}
    \theta_2[f_n] - \theta_2[g_n]
    =
    -
    \mathcal{B}[h_n,g_n]
    -
    \mathcal{B}[g_n,h_n].
\end{align}
We decompose the first term of the right hand side of the above as 
\begin{align}
    \mathcal{B}[h_n,g_n]
    &=
    (-\Delta)^{-\alpha}[\partial_{x_2}(-\Delta)^{-1/2}\theta_1[h_n]\partial_{x_1}\theta_1[g_n]]\\
    &\quad
    -
    (-\Delta)^{-\alpha}[\partial_{x_1}(-\Delta)^{-1/2}\theta_1[h_n]\partial_{x_2}\theta_1[g_n]]\\
    &=:
    \mathcal{B}_n^{(1)} - \mathcal{B}_n^{(2)}.
\end{align}
We further decompose $\mathcal{B}^{(1)}$ as
\begin{align}
    \mathcal{B}_n^{(1)}
    &=
    2^{-(3-2\alpha)n}
    (-\Delta)^{-\alpha}
    \lp{\partial_{x_2}[\phi(x_1)\phi(x_2)]\partial_{x_1}[\phi(x_1)\phi(x_2)\sin(2^nx_1)]}\\
    &=
    \delta^22^{-(2-4\alpha)n}
    (-\Delta)^{-\alpha}\lp{\sp{\phi(x_1)}^2\phi(x_2)\phi'(x_2)\cos(2^nx_1)}\\
    &\quad
    +
    \delta^22^{-(3-4\alpha)n}
    (-\Delta)^{-\alpha}[\phi(x_1)\phi'(x_1)\phi(x_2)\phi'(x_2)\sin(2^nx_1)]\\
    &=:
    \mathcal{B}_n^{(1,1)}+\mathcal{B}_n^{(1,2)}.
\end{align}
Since $\supp \widehat{\mathcal{B}^{(1,1)}_n}\subset \{ \xi \in \mathbb{R}^2 \ ;\  2^n -4 \leq |\xi| \leq 2^n + 4 \}$, we have
\begin{align}
    \n{\mathcal{B}_n^{(1,1)}}_{\dH^{2-2\alpha}}
    \sim{}&
    2^{(2-2\alpha)n}
    \n{\mathcal{B}_n^{(1,1)}}_{L^2}\\
    \sim{}&
    \delta^2\n{\sp{\phi(x_1)}^2\phi(x_2)\phi'(x_2)\cos(2^nx_1)}_{L^2}\\
    ={}&
    \delta^2\n{\sp{\phi(x_1)}^2\cos(2^nx_1)}_{L^2(\mathbb{R}_{x_1})}
    \n{\phi(x_2)\phi'(x_2)}_{L^2(\mathbb{R}_{x_2})}.
\end{align}
Here, we see that $\n{\phi(x_2)\phi'(x_2)}_{L^2(\mathbb{R}_{x_2})}=(1/2)\n{\phi^2}_{\dH^1(\mathbb{R}_{x_2})}$ and
\begin{align}
    \n{\sp{\phi(x_1)}^2\cos(2^nx_1)}_{L^2(\mathbb{R}_{x_1})}^2
    ={}&
    \int_{\mathbb{R}}
    \sp{\phi(x_1)}^4\cos^2(2^nx_1)dx_1\\
    ={}&
    \frac{1}{2}
    \n{\phi}_{L^4(\mathbb{R})}^4
    +
    \frac{1}{2}
    \int_{\mathbb{R}}
    \sp{\phi(x_1)}^4\cos(2^{n+1}x_1)dx_1\\
    \to{}& 
    \frac{1}{2}
    \n{\phi}_{L^4(\mathbb{R})}^4
\end{align}
as $n \to \infty$, where we have used the Riemann--Lebesgue lemma for the last limit.
Thus we have 
\begin{align}
    \liminf_{n \to \infty}
    \n{\mathcal{B}_n^{(1,1)}}_{\dH^{2-2\alpha}}
    \geq c_0\delta^2
\end{align}
for some positive constant $c_0=c_0(\alpha)$.
For the estimate of $\mathcal{B}_n^{(1,2)}$, we have by $\supp \widehat{\mathcal{B}^{(1,2)}_n}\subset \{ \xi \in \mathbb{R}^2 \ ;\  2^n -4 \leq |\xi| \leq 2^n + 4 \}$ that
\begin{align}
    \n{\mathcal{B}_n^{(1,2)}}_{\dH^{2-2\alpha}}
    \leq{}&
    C2^{(2-2\alpha)n}
    \n{\mathcal{B}_n^{(1,2)}}_{L^2}\\
    \leq{}&
    C\delta^22^{-n}
    \n{\phi(x_1)\phi'(x_1)\phi(x_2)\phi'(x_2)\sin(2^nx_1)}_{L^2}\\
    \leq{}&
    C\delta^22^{-n}\n{\phi^2}_{\dH^1(\mathbb{R})}^2.
\end{align}
For the estimate of $\mathcal{B}_n^{(2)}$, we see that 
\begin{align}
    \n{\mathcal{B}_n^{(2)}}_{\dH^{2-2\alpha}}
    \leq{}&
    C
    \n{\partial_{x_1}(-\Delta)^{-1/2}\theta_1[g_n]\partial_{x_2}\theta_1[h_n]}_{\dH^{2-4\alpha}}\\
    \leq{}&
    C\delta^22^{-(3-4\alpha)n}
    \n{\phi(x_1)\phi'(x_1)\phi(x_2)\phi'(x_2)\sin(2^nx_1)}_{\dH^{2-4\alpha}}\\
    \leq{}&
    C\delta^22^{-n}
    \n{\phi(x_1)\phi'(x_1)\phi(x_2)\phi'(x_2)\sin(2^nx_1)}_{L^2}\\
    \leq{}&
    C\delta^22^{-n}\n{\phi^2}_{\dH^1(\mathbb{R})}^2.
\end{align}
For the estimate of $\mathcal{B}[g_n,h_n]$, we have 
\begin{align}
    \n{\mathcal{B}[g_n,h_n]}_{\dH^{2-2\alpha}}
    \leq{}&
    C
    \n{\nabla^{\perp}(-\Delta)^{-1/2}g_n \cdot \nabla h_n}_{\dH^{2-4\alpha}}\\
    \leq{}&
    C
    \n{\theta_1[g_n]}_{\dH^{\alpha}}
    \n{\theta_1[h_n]}_{\dH^{4-5\alpha}}
    +
    C
    \n{\theta_1[g_n]}_{\dH^{2-2\alpha}}
    \n{\theta_1[h_n]}_{\dH^{2-2\alpha}}\\
    \leq{}&
    C\delta^2 2^{-(1-2\alpha)n}.
\end{align}
Hence, we obtain 
\begin{align}
    &\liminf_{n\to \infty}
    \n{\theta_2[f_n]-\theta_2[g_n]}_{\dH^{2-2\alpha}}
    \\
    &\quad\geq
    \liminf_{n\to \infty}
    \sp{
    \n{\mathcal{B}_n^{(1,1)}}_{\dH^{2-2\alpha}}
    -
    \n{\mathcal{B}_n^{(1,2)}}_{\dH^{2-2\alpha}}
    -
    \n{\mathcal{B}_n^{(2)}}_{\dH^{2-2\alpha}}
    -
    \n{\mathcal{B}[g_n,h_n]}_{\dH^{2-2\alpha}}
    }\\
    &\quad\geq c_0\delta^2.
\end{align}
\subsection*{Step 3. Estimates for the remainder terms}
Let us define the remainder terms as
\begin{align}
    &\Theta[f_n]:=\theta[f_n]-\theta_1[f_n]-\theta_2[f_n],\\ 
    &\Theta[g_n]:=\theta[g_n]-\theta_1[g_n]-\theta_2[g_n].
\end{align}
Then, $\Theta[g_n]$ should solve
\begin{align}
    (-\Delta)^{\alpha}\Theta[g_n]
    &+
    v_1[g_n] \cdot \nabla \theta_2[g_n]
    +
    v_2[g_n] \cdot \nabla \theta_1[g_n]
    +
    v_2[g_n] \cdot \nabla \theta_2[g_n]\\
    &
    +
    V[g_n] \cdot \nabla \theta_1[g_n]
    +
    V[g_n] \cdot \nabla \theta_2[g_n]\\
    &
    +
    v_1[g_n] \cdot \nabla \Theta[g_n]
    +
    v_2[g_n] \cdot \nabla \Theta[g_n]
    +
    V[g_n] \cdot \nabla \Theta[g_n]=0,
\end{align}
where we have defined $v_m[g_n]:=\nabla^{\perp}(-\Delta)^{-1/2}\theta_m[g_n]$ with $m=1,2$ and $V[g_n]:=\nabla^{\perp}(-\Delta)^{-1/2}\Theta[g_n]$.
Then from the energy argument as in the proof of Theorem \ref{thm:1}, we have
\begin{align}
    \n{\Theta[g_n]}_{\dH^{\alpha}}
    \leq{}&
    \n{v_1[g_n] \cdot \nabla \theta_2[g_n]
    +
    v_2[g_n] \cdot \nabla \theta_1[g_n]
    +
    v_2[g_n] \cdot \nabla \theta_2[g_n]}_{\dH^{-\alpha}}\\
    &
    +
    \n{V[g_n] \cdot \nabla \theta_1[g_n]
    +
    V[g_n] \cdot \nabla \theta_2[g_n]}_{\dH^{-\alpha}}\\
    \leq{}&
    C\n{\theta_1[g_n]}_{\dH^{\alpha}} \n{\theta_2[g_n]}_{\dH^{2-2\alpha}}
    \\
    &
    +
    C\n{\theta_2[g_n]}_{\dH^{\alpha}} \n{\theta_1[g_n]}_{\dH^{2-2\alpha}}
    \\
    &
    +
    C\n{\theta_2[g_n]}_{\dH^{\alpha}} \n{\theta_2[g_n]}_{\dH^{2-2\alpha}}\\
    &
    +
    C\n{\Theta[g_n]}_{\dH^{\alpha}} \n{\theta_1[g_n]}_{\dH^{2-2\alpha}}
    \\
    &
    +
    C\n{\Theta[g_n]}_{\dH^{\alpha}} \n{\theta_2[g_n]}_{\dH^{2-2\alpha}}.
\end{align}
Here, we see that 
\begin{align}
    &
    \n{\theta_1[g_n]}_{\dH^{\alpha}}
    \leq C2^{-(2-3\alpha)n}\delta, \\
    &
    \n{\theta_1[g_n]}_{\dH^{2-2\alpha}}
    \leq C\delta, \\
    &
    \n{\theta_1[g_n]}_{\dH^{4-5\alpha}}
    \leq C2^{(2-3\alpha)n}\delta, \\
    &
    \n{\theta_1[g_n]}_{\dH^{6-7\alpha}}
    \leq C2^{2(2-3\alpha)n}\delta, \\
\end{align}
and
\begin{align}
    &
    \n{\theta_2[g_n]}_{\dH^{\alpha}}
    \leq
    C
    \n{\theta_1[g_n]}_{\dH^{\alpha}} \n{\theta_1[g_n]}_{\dH^{2-2\alpha}}
    \leq
    C2^{-(2-3\alpha)n}\delta^2,\\
    &
    \n{\theta_2[g_n]}_{\dH^{2-2\alpha}}
    \leq
    C\n{\theta_1[g_n]}_{\dH^{2-2\alpha}}^2
    +
    C
    \n{\theta_1[g_n]}_{\dH^{\alpha}} \n{\theta_1[g_n]}_{\dH^{4-5\alpha}}
    \leq
    C\delta^2,\\
    &
    \begin{aligned}
    \n{\theta_2[g_n]}_{\dH^{4-5\alpha}}
    \leq{}&
    C\n{\theta_1[g_n]}_{\dH^{2-2\alpha}}\n{\theta_1[g_n]}_{\dH^{4-5\alpha}}
    +
    C
    \n{\theta_1[g_n]}_{\dH^{\alpha}} \n{\theta_1[g_n]}_{\dH^{6-7\alpha}}
    \\
    \leq{}&
    C2^{(2-3\alpha)n}\delta^2.
    \end{aligned}
\end{align}
Thus, choosing $\delta$ sufficiently small, we have 
\begin{align}
    \n{\Theta[g_n]}_{\dH^{\alpha}}
    \leq
    C2^{-(2-3\alpha)n}\delta^3.
\end{align}
For the energy estimates in the $\dH^{2-2\alpha}(\mathbb{R}^2)$-norm, we have
\begin{align}
    \n{\Theta[g_n]}_{\dH^{2-2\alpha}}
    \leq{}&
    \n{v_1[g_n] \cdot \nabla \theta_2[g_n]
    +
    v_2[g_n] \cdot \nabla \theta_1[g_n]
    +
    v_2[g_n] \cdot \nabla \theta_2[g_n]}_{\dH^{2-4\alpha}}\\
    &
    +
    \n{V[g_n] \cdot \nabla \theta_1[g_n]
    +
    V[g_n] \cdot \nabla \theta_2[g_n]}_{\dH^{2-4\alpha}}\\
    &
    +
    \n{V[g_n] \cdot \nabla \theta_1[g_n]
    +
    V[g_n] \cdot \nabla \theta_2[g_n]}_{\dH^{2-4\alpha}}\\
    &
    +
    \sum_{m=1}^2
    \n{\lp{v_{m}[g_n]\cdot \nabla, (-\Delta)^{\frac{2-3\alpha}{2}}} \theta_N}_{\dH^{-\alpha}}\\
    &
    +
    \n{\lp{V[g_n]\cdot \nabla, (-\Delta)^{\frac{2-3\alpha}{2}}} \theta_N}_{\dH^{-\alpha}}
    \\
    \leq{}&
    C\n{\theta_1[g_n]}_{\dH^{\alpha}} \n{\theta_2[g_n]}_{\dH^{4-5\alpha}}
    \\
    &
    +
    C\n{\theta_2[g_n]}_{\dH^{\alpha}} \n{\theta_1[g_n]}_{\dH^{4-5\alpha}}
    \\
    &
    +
    C\n{\theta_2[g_n]}_{\dH^{\alpha}} \n{\theta_2[g_n]}_{\dH^{4-5\alpha}}\\
    &
    +
    C\n{\Theta[g_n]}_{\dH^{\alpha}} \n{\theta_1[g_n]}_{\dH^{4-5\alpha}}
    \\
    &
    +
    C\n{\Theta[g_n]}_{\dH^{\alpha}} \n{\theta_2[g_n]}_{\dH^{4-5\alpha}}\\
    &
    +
    C
    \n{\theta_1[g_n]}_{\dH^{2-\alpha}}
    \n{\Theta[g_n]}_{\dH^{2-2\alpha}}\\
    &
    +
    C
    \n{\theta_2[g_n]}_{\dH^{2-\alpha}}
    \n{\Theta[g_n]}_{\dH^{2-2\alpha}}\\
    &
    +
    C
    \n{\Theta[g_n]}_{\dH^{2-2\alpha}}
    \n{\Theta[g_n]}_{\dH^{2-2\alpha}},
\end{align}
which and all above estimates in this step imply 
\begin{align}
    \n{\Theta[g_n]}_{\dH^{2-2\alpha}}
    \leq
    C\delta^3 + C\delta\n{\Theta[g_n]}_{\dH^{2-2\alpha}}+C\sp{\n{\Theta[g_n]}_{\dH^{2-2\alpha}}}^2.
\end{align}
Since $\n{\Theta[g_n]}_{\dH^{2-2\alpha}} \leq \n{\theta_1[g_n]}_{\dH^{2-2\alpha}}+\n{\theta_2[g_n]}_{\dH^{2-2\alpha}}+\n{\theta[g_n]}_{\dH^{2-2\alpha}}\leq C\delta$, we solve the above quadratic inequality to obtain 
\begin{align}
    \n{\Theta[g_n]}_{\dH^{2-2\alpha}} \leq C\delta^3.
\end{align}
By the similar manner, we have 
\begin{align}
    \n{\Theta[f_n]}_{\dH^{2-2\alpha}} \leq C\delta^3.
\end{align}
Hence, it holds
\begin{align}
    &
    \n{\theta[f_n]-\theta[g_n]}_{\dH^{2-2\alpha}}\\
    &\quad
    \geq
    \n{\theta_2[f_n]-\theta_2[g_n]}_{\dH^{2-2\alpha}}
    -
    \n{\theta_1[f_n]-\theta_1[g_n]}_{\dH^{2-2\alpha}}\\
    &\qquad
    -
    \n{\Theta[f_n]}_{\dH^{2-2\alpha}}
    -
    \n{\Theta[g_n]}_{\dH^{2-2\alpha}}\\
    &\quad 
    \geq
    \n{\theta_2[f_n]-\theta_2[g_n]}_{\dH^{2-2\alpha}}
    -
    C
    \n{f_n-g_n}_{\dH^{2-4\alpha}}
    -
    C\delta^3,
\end{align}
which and the results in the previous step yield
\begin{align}
    \liminf_{n \to \infty}
    \n{\theta[f_n]-\theta[g_n]}_{\dH^{2-2\alpha}}
    \geq c_0\delta^2-C\delta^3.
\end{align}
Thus, choosing $\delta$ sufficiently small, we complete the proof.
\end{proof}



\noindent
{\bf Data availability.}\\
Data sharing not applicable to this article as no datasets were generated or analysed during the current study.

\noindent
{\bf Conflict of interest.}\\
The author has declared no conflicts of interest.

\noindent
{\bf Acknowledgements.} \\
The author was supported by Grant-in-Aid for Research Activity Start-up, Grant Number JP2319011.
The author would like to express his sincere gratitude to Professor Hideo Kozono for many fruitful advices and continuous encouragement.

\end{document}